\documentclass[12pt]{article}

\usepackage{graphicx}
\usepackage{amsmath,amssymb}
\usepackage{subfig}

\begin{document}

\title{\Large \bf Approximating cube roots of integers, \\
 after Heron's {\em Metrica III.20\/} }

\author{{\sc Trond Steihaug and D. G. Rogers} \\
  {\small \em Institutt for Informatikk, Universitetet i Bergen} \\
  {\small \em PB7803, N5020, Bergen, Norge} \\
  {\small \em {\tt trond.steihaug@ii.uib.no\/} } }

\date{}

\maketitle

\renewcommand{\baselinestretch}{1.05}

\setlength{\parindent}{0pt}
\setlength{\parskip}{1.0ex}

\begin{center}\em \small
For Christian Marinus Taisbak, \\
Institut for Gr{\ae}sk og Latin,
K{\o}benhavns Universitet, 1964--1994, \\
On his eightieth birthday, 17 February, 2014
\end{center}

\begin{quote}\em
Heron did not need any other corroboration than the fact that
 the method works, and that the separate results are easily
 confirmed by multiplication.
\begin{flushright}\em
C.~M.~Taisbak \cite[\S 2]{CMT}
\end{flushright}
\end{quote}

\section{Taisbak's conjecture}

How often, in the happy Chinese idiom, do we search high and
 low for our shoulder pole, only at last to notice it again
 on our shoulder where we left it?
For all that the learned commentator might reassure us that 
 some mathematician of the past could not help but make some
 pertinent observation, just as surely we know, from our own
 experience, that such acuity might escape us for half a
 lifetime, before, all at once, perhaps of a Summer's night,
 the {\o}re drops.
This is, indeed, the story behind Christian Marinus Taisbak's
 conjecture in \cite{CMT}, as divulged in a recent letter
 \cite{CMT2}.
So, we too were set thinking.
We report here on some of our findings.

Heron, in {\em Metrica III.20--22\/}, is concerned with the
 the division of solid figures --- pyramids, cones and {\em
 frustra\/} of cones --- to which end there is a need to
 extract cube roots \cite[II, pp.~340--342]{TLH} (see also
 \cite[p.~430]{TLHM}).
A case in point is the cube root of 100, for which Heron
 obligingly outlines a method of approximation in {\em
 Metrica III.20\/} as follows (adapted from \cite{TLH,OB,CMT},
 noting that the addition in \cite[p.~103, fn.~1]{CMT} appears
 earlier in \cite[p.~69]{OB}; cf.\
 \cite[p.~191, fn.~124]{WRK86}):
\begin{quote}\small \em
Take the cube numbers nearest 100 both above and below, namely
 125 and 64. \\
Then, $125-100=25$ and $100-64=36$. \\
Multiply 25 by 4 and 36 by 5 to get 100 and 180; and then add
 to get 280. \\
Divide 180 by 280, giving 9/14. Add this to the side of the
 smaller cube; this gives $4\frac{9}{14}$ as the cube root of
 100 as nearly as possible.
\end{quote}

It seems short, unobjectionable work to turn this descriptive
 algorithm into a general formula for approximating the cube
 root of some given integer $N$.
We first locate $N$ among the cubes of the integers:
$$ m^3 < N < (m+1)^3. $$
Writing $d_1=N-m^3$ and $d_2=(m+1)^3-N$, Heron would then have
 us approximate the cube root of $N$ by
\begin{equation}
 m+\frac{(m+1)d_1}{(m+1)d_1+md_2}.
\end{equation}

The text of {\em Metrica\/} as we have it today only came
 to light in the mid-1890s, with a scholarly edition
 \cite{HS} published in 1903.
How little was known for sure about {\em Metrica\/} in the
 years immediately prior to this is suggested by \cite{JG}.
Fragments were known by quotation in other sources and
 Eutocius, in a commentary on the works of Archimedes,
 reports that Heron used the same methods for square and
 cube roots as Archimedes.
But clearly this does not have the same {\em cachet\/} as a
 text --- and we still lack anything by Archimedes on finding
 cube roots.
Gustave Wertheim (1843--1902) proposed (1) in 1899 in \cite{GW},
 to be followed a few years latter by Gustaf Hjalmar
 Enestr{\" o}m (1852--1923) in \cite{GHE2} with an exact (if
 tautological) expression, given below in \S 5.2 as (22), for
 the cube root of $N$ from which (1) follows on discarding
 cubes of positive terms less than unity.
(Besides work in mathematics and statistics, Enestr{\" o}m had
 interests in the history of mathematics, as seen, for example,
 in his note \cite{GHE} on rules of convergence in the 1700s:
 he is perhaps best remembered today for introducing the
 {\em Enestr{\" o}m index\/} to help identify the writings of
 Leonhard Euler (1707--1783); but, while there seems to be
 little written about him in English, the very first volume
 of {\em Nordisk Matematisk Tidskrift\/} carried a centenary
 profile \cite{KGH}.)
 
To be sure, other formulae might fit Heron's numerical instance
 in {\em Metrica III.20\/}: a nod is made to one in
 \cite[pp.~137--138]{AH}:
\begin{equation}
   m+\frac{d_1\sqrt{d_2}}{N+d_1\sqrt{d_2}}.
\end{equation}
At first sight, this gesture might seem {\em pro forma\/}, as
 it is conceded straightaway that (2), when compared with (1),
 is both less easy to justify and not so accurate for other
 values of $N$.
But the record has not always been so clear-cut and it is (2),
 not (1), that we find on looking back to
 \cite[pp.~62--63]{IT}, where reference is made to an article
 \cite{MC} by Ernst Ludwig Wilhelm Maximilian Curtze
 (1837--1903) of 1897, along with \cite{GW,GHE2}.
Both Curtze's tentative contribution (2) and another, similar
 formula,
\begin{equation}
  m+\frac{(m+1)d_1}{N+(m+1)d_1},
\end{equation}
 had, in fact, been compared adversely for accuracy with (1)
 in 1920 by Josiah Gilbart Smyly (1867--1948) in \cite{JGS};
 Smyly attributes to George Randolph Webb (1877--1929; Fellow,
 Trinity College, Dublin) an estimate that the error in
 (1) is of the order of $1/m^2$ (see further \S 5.2,
 especially (33)).
For the record, we might note here that Smyly alludes to the
 work of Curtze, but not that of Wertheim or Enestr{\" o}m;
 on the other hand, Heath \cite{TLH} cites them, but not
 Curtze or Smyly (truely the vagaries of citation are not
 easily explicable: in \cite[p.~256, fn.~2]{DES}, we find
 Smyly footnoted as {\em correcting\/} Curtze, only for (3),
 rather than (1), to be printed).

It is also worth observing that the effect of emendations is
 to move our understanding of the received text in favour of
 the most accurate candidate, namely (1). 
As it happens, in {\em Metrica III.22\/}, Heron needs the
 estimate of another cube root, that of 97050 according to
 \cite{HS}, but in fact of $97804\frac{4}{5}$, as pointed
 out in \cite[pp.~338--340]{EMB}. 
The approximation taken is 46, which cubes to 97336, so is
 not too far off either way, suggesting that Heron did not
 allow himself to be blinded by science.

If the consensus on (1) is by now reasonably settled, there
 remains the question of how Heron might have come upon (1),
 as well as the somewhat different question of how (1) might be
 justified.
A formal derivation of (1) might well fail to satisfy those who
 want some heuristic insight into the approximation; and
 Enestr{\" o}m may have lost sight of the simplicity of his
 identity (26) in the manner he derives it (see further \S 5.2).
Taisbak strikes out on his own account in \cite{CMT} from the
 constancy of the third difference of the sequence of cubes of
 integers and builds up to the observation that the gradient of
 the chord between $m-1$ and $m$ is to the gradient of the
 chord between $m$ and $m+1$ approximately as $m-1:m+1$.
In effect, Taisbak sums up his thinking with a question
 \cite[\S 3]{CMT}: ``{\em Did the Ancients know and use
 sequences of differences?\/}''

As far as Taisbak's mathematics goes, a rather similar argument
 was advanced some thirty years ago by Henry Graham Flegg
 (1924-- ) in a book \cite[p.~137]{HGF} (pleasingly enough it
 was reissued in 2013).
Others have been here, too: Oskar Becker (1889--1964) in
 \cite[pp.~69--71]{OB} in 1957; Evert Marie Bruins (1909-1990)
 in \cite[p.~336]{EMB} in 1964; Wilbur Richard Knorr (1945--1997)
 in \cite[pp.~191--194]{WRK86} in 1986.
It has also been noted how (1) can be adapted for iterative use,
 although the accuracy of (1), as remarked on by Smyly, coupled
 with the opportunity for rescaling it provides, might make
 iteration otiose (cf.\ \S 5.4).
But, in fairness, it might be remarked that the main difference
 between these writers and Enestr{\" o}m is that their
 approximative sleight of hand takes good care to wipe away
 small terms as they go, rather than in one fell swoop at the
 end (we return to these comparisons in \S 5.2). 

Our concerns are rather different.
For a start, might there be more to discern in the numerical
 instance Heron presents in {\em Metrica\/}?
This prompts two further questions.
Why is no comparison made with the more straightforward cube
 root bounds (as in (6) and (7)) analogous to those (as in
 (14) and (15)) seemingly in common use by Archimedes, Heron
 and others for square roots?
And, why do we not hear anything like (1) in regard to square
 roots?
Then, again, might there not be more to say about (1) itself?

Our concerns in these regards are mathematical, not historical.
Perforce, we respect Taisbak's stricture, as endorsed by
 Unguru \cite{SU}, that we adopt as our epigraph.
Truly, the proof of the pudding is in the eating; and if,
 perhaps like Eutocius commenting on the works of Archimedes,
 you have nothing more imaginative to offer, arithmetical
 confirmation remains a safe recourse, if not always a sure one
 (cf.\ \cite[pp.~522, 540]{WRK89}).
But we suspect that, if anything, others before us may have been
 too abashed to descend our level of {\em naivet{\' e}\/}.
Our excuse, if one is needed, is that, even at this level,
 there is still much with which to be usefully engaged.
 
\section{Heron's example}

The difference between successive cubes is
\begin{equation}
  (m+1)^3-m^3=3m^2+3m+1.
\end{equation}

\begin{figure}[htp]

\begin{center}

\includegraphics[width=8cm]{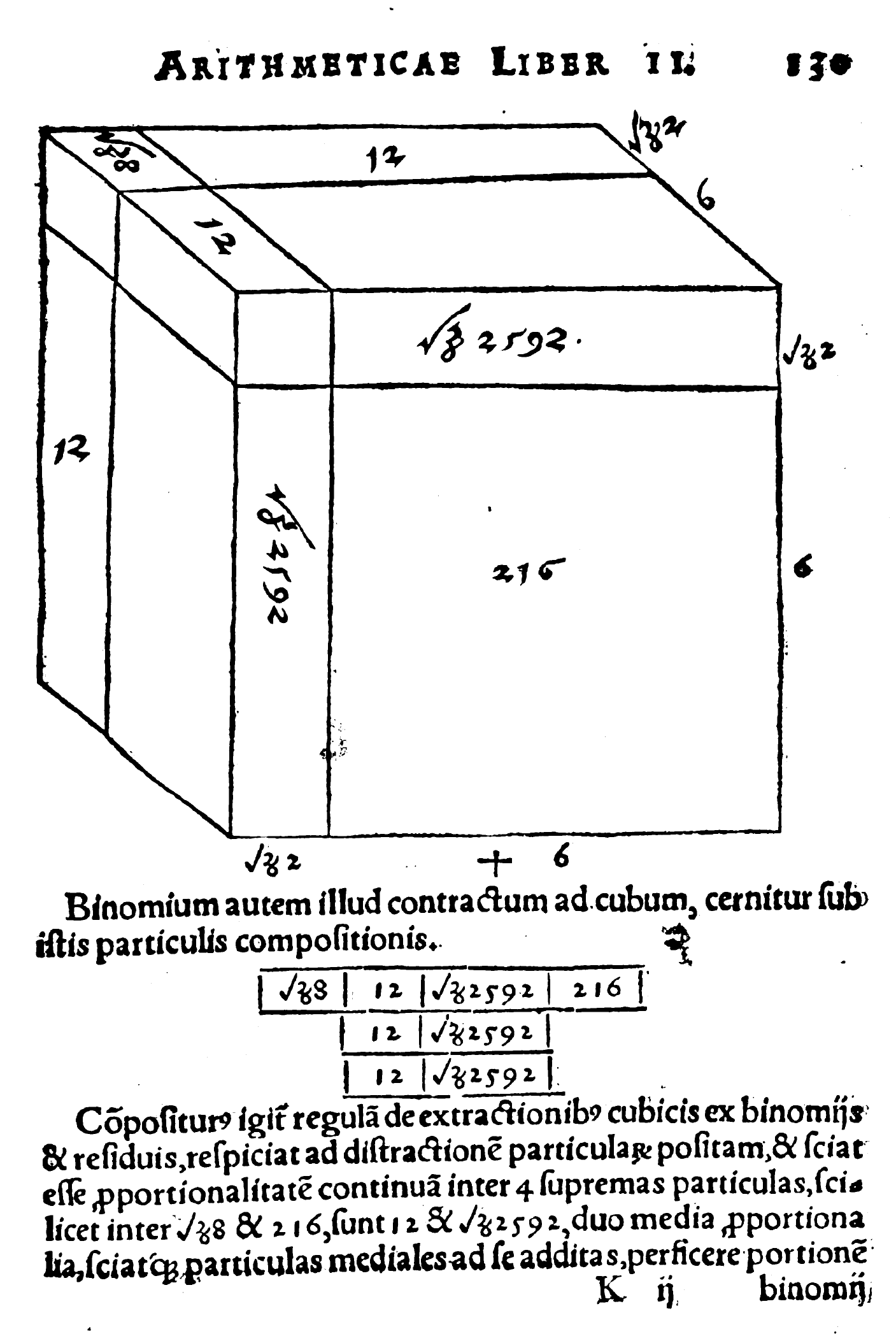}
\caption{Picture of a cubed binomial from {\em Arithmetica Integra\/} (1544)}

\end{center}
\end{figure}

More generally, we may picture the difference between cubes
 by cutting up the larger cube into smaller cube with various
 other slabs and blocks, a three-dimensional analogue of the
 pictures we might draw for the difference of two squares,
 perhaps as an {\em aide m{\' e}moire\/} to our reading of
 Euclid's {\em Elements II\/} (one traditional mode of
 visualising the cube of a binomial expression is shown in
 Fig.~1; an alternative dissection appears in Fig.~3 in
 conjunction with (16)).

Thus, as $d_1=N-m^3$ and $d_2=(m+1)^3-N$ sum to this
 difference, we can ensure some cancellation in working with
 (1) if we arrange to take $d_1$ to be $k(m+1)+1$ for some
 $k$ with $0\leq k\leq 3m$.
Perhaps Heron had something of this in mind in taking an
 example in which $d_1=(2m-1)(m+1)+1=m(2m+1)$ and
 $d_2=(m+1)^2$ for $m=4$.
At all events, generalising Heron's example in this way, we
 obtain from (1) a bound on the cube root of
 $N=m^3+m(2m+1)=m(m+1)^3-(m+1)^2$:
\begin{equation}
   m+\frac{2m+1}{3m+2}=m+1-\frac{m+1}{3m+2}.
\end{equation}
It is a simple matter of verification to check that this is
 an {\em upper\/} bound.

But not only is this pleasing in itself, the form of these
 expressions suggests --- invites? --- a comparison with the
 upper bounds obtained more straightforwardly from binomial
 expressions analogous to those familiar for square roots
 (as in (14) and (15)), of which Gerolamo Cardono (1501--1576)
 made celebrated use in {\em Practica Arithmetice\/} (1539)
 \cite[\S 2.4]{HLH} (but cf.\ also (20)).
Thus, for $N=m^3+d_1$, the cube root is bounded above by
\begin{equation}
   m+\frac{d_1}{3m^2},
\end{equation}
 while for $N=(m+1)^3-d_2$, the cube root is bounded above by
\begin{equation}
   m+1-\frac{d_2}{3(m+1)^2}.
\end{equation}

So, in generalising Heron's example, we have hit on a case
 where the upper bounds in (6) and (7) also come out rather
 neatly:
$$ m+\frac{2m+1}{3m};\quad m+1-\frac{1}{3}. $$
Of course, the former is not so good as the latter, reflecting
 the closer proximity of {\em this\/} $N$ to $(m+1)^3$ than
 to $m^3$.
Rather more strikingly neither of these bounds is as good as
 that in (5) obtained from (1); indeed,
$$ \frac{2m+1}{3m+2} < \frac{2}{3} < \frac{2m+1}{3m}. $$
It is possible to squeeze (6) further by increasing the
 denominator in the fraction, and some writers in Arabic in the
 early 1000s worked with $3m^2+1$ in place of $3m^2$ (cf.\
 \cite[\S 3.2]{HLH}.
But this still does not give an improvement over (5).

Whether or not Heron may have indulged himself in such
 exercises, a few numerical instances like this would surely
 convey to any impressionable mind that (1) cannot be
 completely without merit.
Trouble might spring more from the opposite corner, not to run
 away with too favourable an endorsement based {\em only\/} on
 evidence of this sort.
However, as we show in \S 5.3, an approximate construction of
 two mean proportionals examined by Pappus early in {\em
 Synagogue III} allows us to improve on (5), indicating that
 it is by no means the best the Greeks could have done, had
 they put their minds to it. 

\section{Square roots}

\subsection{Elementary theory of proportions\/}

When we look at the formulation of (1), it would seem that it is
 a recipe we could write down for other functions besides cubes
 and cube roots; and, if for cubes and cube roots, why not before
 that for squares and square roots?
In fact, we might recognize (1) in the setting of the elementary
 theory of proportions that was well-articulated by the Greeks.
For, given $a/b>c/d>0$, an early result in that theory gives
$$ \frac{c}{d} < \frac{a+c}{b+d} < \frac{a}{b}, $$
 and, more generally, for weights $w_1$ and $w_2$,
\begin{equation}
   \frac{c}{d} < \frac{aw_1+cw_2}{bw_1+dw_2} < \frac{a}{b}.
\end{equation}
In particular (cf.\ (11), (13), (22) and (26))
\begin{equation}
   m=\frac{m^2}{m}<\frac{(m+1)^2w_1+m^2w_1}{(m+1)w_1+mw_2}
    <\frac{(m+1)^2}{m+1}=m+1,
\end{equation}
 where the central expression can then be rewritten as (cf.\ (1))
$$ \frac{(m+1)^2w_1+m^2w_1}{(m+1)w_1+mw_2}
  =m+\frac{(m+1)w_1}{(m+1)w_1+mw_2}. $$
This is pudding that anyone can eat, but it might not always satisfy
 Winston Churchill's demand that pudding have a {\em theme\/}.
For, how to explain the choice of weights for different functions?

\subsection{Curves and chords\/}

\begin{figure}[htp]

\begin{center}

\includegraphics[width=8cm]{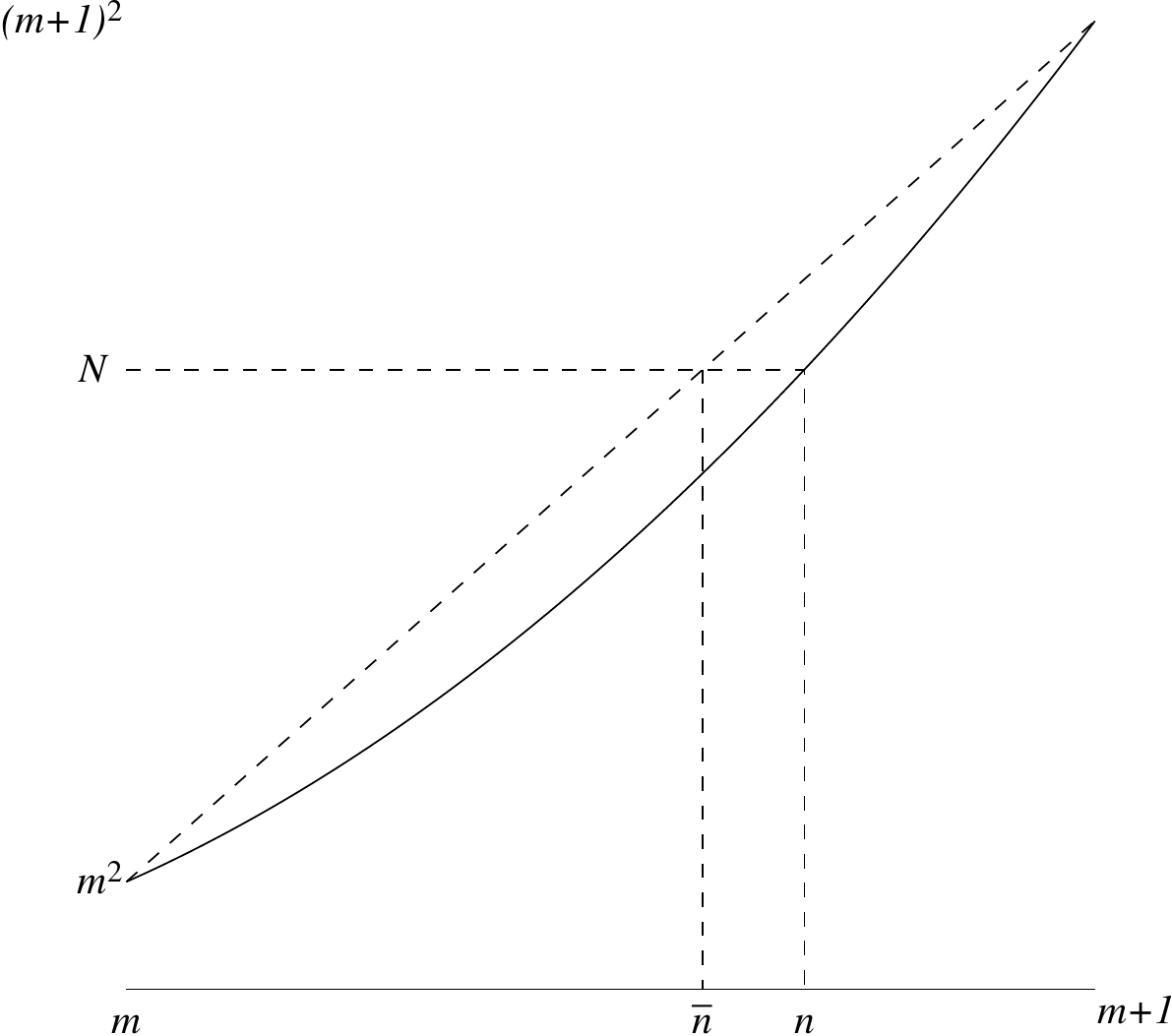}
\caption{Approximating square roots from below}

\end{center}
\end{figure}

For any increasingly increasing function, such as squaring or
 cubing, chords lie above the curve, so a particular height $N$
 will be encountered on the chord before it is encountered on the
 curve, giving a simple means of finding a lower bound on the
 ordinate for which $N$ is attained, after the manner of solution
 traditionally known as ``{\em double false position\/}'' (a brief
 introduction to the history of which is recently to hand in
 \cite{H}).
Let us illustrate the thinking here rather naively in the case of
 squares.
So, suppose now that we are given $N$, with
$$ m^2 < N < (m+1)^2, $$
 and we are interested in the square root $n=\sqrt{N}$.
Then we expect that the gradient of the chord between $m$ and $n$,
 that is, $d_1/(n-m)$, to be less than the gradient of the chord
 going on from $n$ to $m+1$, that is, $d_2/(n+1-m)$, where for our
 present purposes in this section we write $d_1=n^2-m^2$ and
 $d_2=(m+1)^2-n^2$ in analogy with the notation for (1).
But, if
\begin{equation}
  \frac{d_1}{n-m} < \frac{d_2}{m+1-n},
\end{equation}
 then it follows that, for $0\leq d_1\leq 2m+1$,
\begin{equation}
  n > \frac{(m+1)d_1+md_2}{d_1+d_2}=m+\frac{d_1}{2m+1}.
\end{equation}
Equality would hold here if the two gradients were equal, in which
 case the common value would be the gradient of the chord from $m$
 to $m+1$, confirming that this lower bound on $n$ is the ordinate
 $\bar{n}$ at which $N$ is attained on this chord (as in Fig.~2).

Of course, in this case, $d_1$ and $d_2$ are just differences of
 squares,
$$ d_1=n^2-m^2=(n-m)(n+m);\quad d_2=(m+1)^2-n^2=(m+1-n)(m+1+n), $$ 
so
$$ \frac{d_1}{n-m}=n+m; \quad \frac{d_2}{m+1-n}=m+1+n. $$
Hence, (10) holds trivially:
$$ n+m < m+1+n. $$

But, looking at this last inequality, we see that it is readily
 reversed by judicious counterpoised weighting, mutiplying the
 left-hand side by $m+1$ and the right hand side by $m$:
$$ (m+1)(n+m) > m(m+1+n). $$
So, in addition to (10), we also have
\begin{equation}
   \frac{(m+1)d_1}{n-m} > \frac{md_2}{m+1-n}
\end{equation}
 from which we deduce in turn the upper bound
\begin{equation}
   n < \frac{(m+1)^2d_1+m^2d_2}{(m+1)d_1+md_2}=m+\frac{(m+1)d_1}{(m+1)d_1+md_2},
\end{equation}
 thereby providing easy confirmation that the analogue of (1)
 for square roots.
But the algebra here is such that conversely, if a upper bound of
 the form (13) holds, then the weighted gradients stand as in (12),
 a point to bear in mind when considering (1).

\subsection{Square root bounds\/}

However, the sad fact of the matter is that (13) is not much help
 because we already do better with one or other of the standard
 upper bounds for square roots obtained from binomial expressions
 that complement the lower bound (11); the implicit use of all
 the bounds (11), (14) and (15) in antiquity is examined {\em in
 extenso\/} in \cite[pp.~lxxvii--xcix]{TLHA} (cf.\
 \cite[pp.~53--57]{JG}).
We recall that, for $N=m^2+d_1$,
\begin{equation}
 n=\sqrt{N} < m+\frac{d_1}{2m}
\end{equation}
 while, for $N=(m+1)^2-d_2$,
\begin{equation}
 n=\sqrt{N} < m+1-\frac{d_2}{2(m+1)}.
\end{equation}
We work with (14) for $0 < d_1\leq m$, switching to (15) for
 $0 < d_2\leq m+1$.

Notice that (14) and (15) also follow from the iterative scheme that
 Heron sketches by example for $N=720$ in {\em Metrica I.8\/}:
$$ m_1=\frac{1}{2}\left(\frac{N}{m_0}+m_0\right), $$
 with $m_0=m$ for (14) and $m_0=m+1$ for (15).
Whether Heron recognised (15) explicitly depends in large part on
 what inference can be drawn from the way fractions are recorded
 (cf.\ \cite[II, p.~326]{TLH}).
There are other puzzles in relation to Heronian iteration.
For instance, samplings in \cite[p.~lxxxii]{TLHA} and \cite[p.~6]{WCE}
 of estimates used by Heron for square roots includes that for
 $\sqrt{75}$ as $8\frac{11}{16}$ (cf.\ (14)), rather than
 $8\frac{2}{3}$ (cf.\ (15); and see further \cite[pp.~10--11]{EMB}),
 which is simpler, as well as more accurate; and a further example is
 raised in \S 5.4.  

Now, in these ranges for $d_1$ and $d_2$ for (14) and (15),
$$ (m+1)d_1+md_2\leq 2m(m+1), $$
 with equality if and only if $d_1=m$ and $d_2=m+1$.
Hence (13) is only as good as (14) or (15) in the case where $d_1=m$
 and $d_2=m+1$, when all three bounds come out the same, namely
 $m+\frac{1}{2}$ (but see \S 5.4 for a reprieve of sorts for (9)). 
This points up the altered situation for cube roots, where the
 evidence of the previous section shows that (1) does better
 than (4) and (5), at least in a family of instances generalizing
 Heron's example in {\em Metrica III.20\/}.
Clearly, we need to examine how the arguments leading to (11) and
 (3) for square roots go over to cube roots, especially as it is
 the innocent use of counterpoised weighting in shifting from (10)
 to (12) that lies at the heart of Taisbak's musings in \cite{CMT}.

\subsection{Mellema's formula for quadratics\/}

But before leaving this discussion of square roots it may be
 instructive in comparison with the derivation of Enestr{\" o}m's
 identity (26) to take a brief look at a formula developed by Elcie
 Edouard Leon Mellema (1544--1622) as a {\em baroque\/} example of
 the method of false position (cf.\ \cite{H}).
Suppose that a function $f(x)$ has a root at $n$ with $a < n < b$,
 then, trivially,
$$ (f(n)-f(a))f(b)=(f(n)-f(b))f(a). $$
However, in the case of a quadratic function where the square has
 been completed, that is, where
$$ f(x)=(x+p)^2-q, $$
 rearranging this equation to make $(n+p)^2$ the subject yields
 Mellema's formula:
$$ (n+p)^2=\frac{(a+p)^2f(b)-(b+p)^2f(a)}{f(b)-f(a)}. $$ 

In contrast with (26), from which (1) follows as an approximation,
 the best that can be said of Mellema's formula is that it is a trick
 on him, if not also on any who might be taken in by it, as it just
 recomputes $q$, which we might suppose would be known more swiftly
 on completing the square in the quadratic.

\section{Cube roots}

So, let us now return to cube roots and our initial supposition
 that we are given $N$, with
$$ m^3 < N < (m+1)^3, $$
 and write
$$ d_1=N-m^3;\quad d_2=(m+1)^3-N. $$

\begin{figure}[htp]

\begin{center}

\includegraphics[width=8cm]{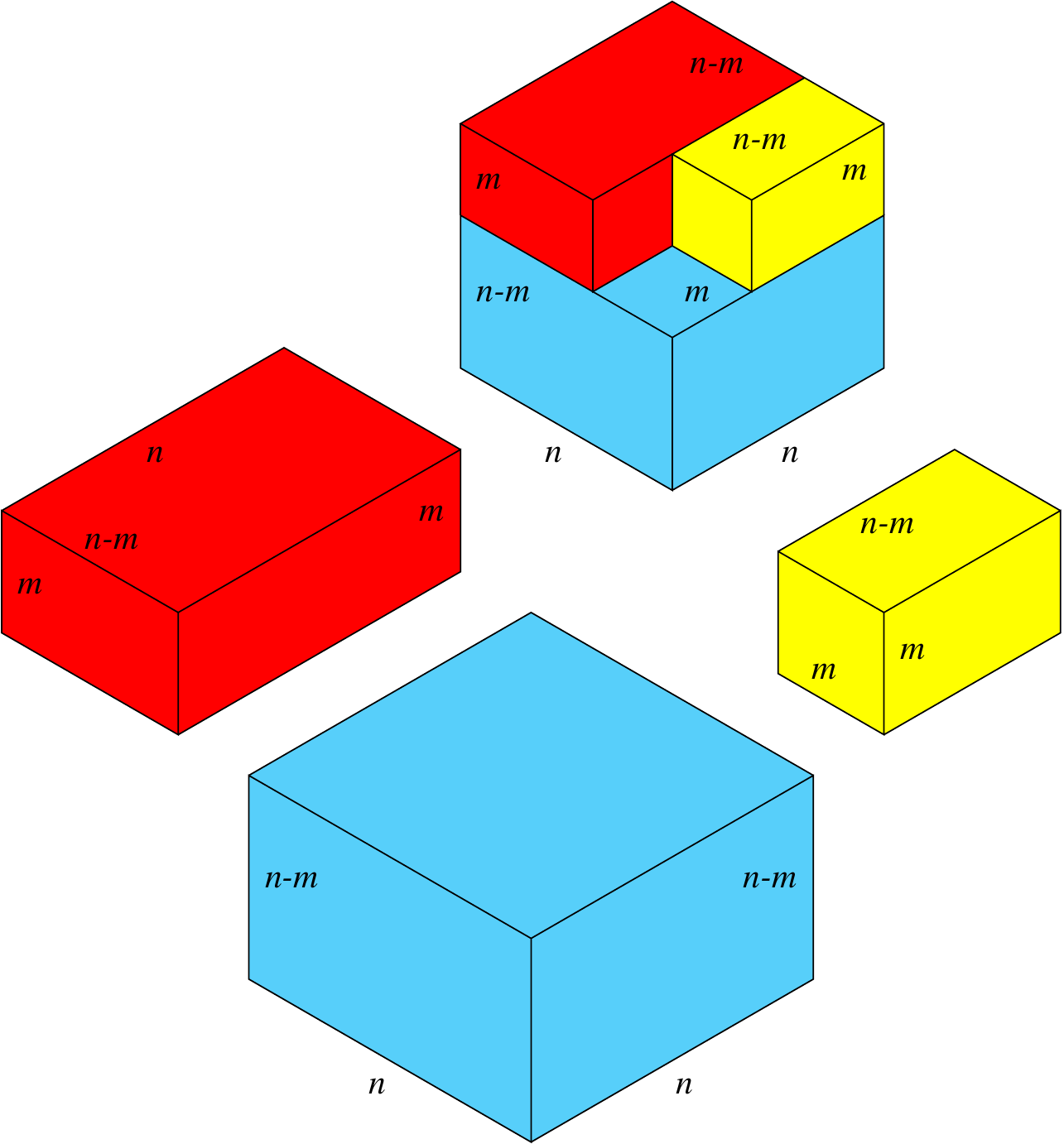}
\caption{Difference of cubes dissected according to (16)}

\end{center}
\end{figure}

If $n$ is the cube root of $N$, so $n^3=N$, then, possibly calling
 to mind Heron's account of {\em frustra\/} of pyramids and cones
 in {\em Metrica II.6, 9\/} (cf. \cite[II, pp.~332--334]{TLH}; that
 the formulae Heron provides were not always used with sufficient
 care is suggested in \cite[pp.~107--108]{JGS2}),
\begin{equation} 
  d_1=n^3-m^3=(n-m)(n^2+nm+m^2),
\end{equation}
 so that
\begin{equation}
  \frac{d_1}{n-m}=n^2+nm+m^2.
\end{equation}
Similarly
\begin{equation}
  \frac{d_2}{m+1-n}=(m+1)^2+(m+1)n+n^2.
\end{equation}
It follows that, on the lines of (10), we have
\begin{equation}
  \frac{d_2}{m+1-n}-\frac{d_1}{n-m}=2m+n+1>0,
\end{equation}
 from which we deduce, in perfect analogy with (11), the lower bound
\begin{equation}
  m_l=\frac{(m+1)d_1+md_2}{d_1+d_2}=m+\frac{d_1}{3m(m+1)+1},
\end{equation}
 and then, iterating the argument, the further refined lower bound
$$  m_l+\frac{N-m_l^3}{3m_l(m+1)+(m+1-m_l)^2}. $$
By way of illustration, in Heron's example with $N=100$, neither the
 the lower bound $m_l=4\frac{36}{61}$ obtained from (20) nor the
 refined one, which involves much heavier computation, are as close
 to the cube root of 100 as Heron's upper bound $4\frac{9}{14}$.
Yet, as a matter of historical record, Leonardo Pisano (Fibonacci;
 1170?--1250?), in {\em Liber Abaci\/} (1202) \cite[\S 2.3]{HLH} and
 again in {\em De Practica Geometrie\/} (1223) \cite[pp.~260--262]{BH},
 approximates cube roots by means of (20), sometimes in sequence with
 its improvement, knowing to ignore the term $(m+1-m_l)^2$ in the
 denominator of the fraction in the latter and even the analogous 1 in
 the denominator of the last fraction in (20) if it suits the
 calculation (the textual problem raised in \cite[p.~92, fn.~7]{HLH}
 as to the use of the improved bound is resolved on cross-reference
 with \cite[p.~262]{BH}).
A version of (20) appears again in use in the 1500s (cf.\
 \cite[p.~255, fn.~4]{DES}; \cite{H}).

So far, so good, although this is entirely as we might expect.
But what about applying Taisbak's hunch on counterpoised
 weightings to (17) and (18) that, as we have seen in the previous
 section, does lead in the case of square roots to the analogue (13)
 of (1)?

Thus, in place of (19), we shall need to consider: 
\begin{equation}
 \frac{(m+1)d_1}{n-m}-\frac{md_2}{m+1-n}=n^2-m(m+1).
\end{equation}
Now, with (21), we see the contingent nature of the expression in
 (1) as a bound on the cube root of $N$.
For, if $N^2>m^3(m+1)^3$, as is certainly the case when
 $N>(m+\frac{1}{2})^3$, then the right-hand side of (21) is positive,
 and, as, in the previous section, it follows that (1) gives an
 upper bound.
On the other hand, if $N^3<m^3(m+1)^3$, (1) will give another lower
 bound along with (20), although one that improves on (20), as it
 is a matter of easy algebra to check that the expression in (1)
 is always larger than its counterpart in (20):
$$ (a^2p+b^2q)(p+q)\geq (ap+bq)^2. $$
In this latter case, let us take by way of illustration $N=85$, so
 $d_1=21$ and $d_2=40$; the two lower bounds then come out as
 $4\frac{21}{61}$, for (20), and $4\frac{21}{53}$, for (1). 

Of course, we can always up the ante by further loading the weights.
Moving up from (21), we find that
$$ \frac{(m+1)^2d_1}{n-m}-\frac{m^2d_2}{m+1-n}=(2m+1)n^2+m(m+1)n>0, $$
 so at least we have the upper bound
\begin{equation}
   n < \frac{(m+1)^3d_1+m^3d_2}{(m+1)^2d_1+m^2d_2},
\end{equation}
 {\em throughout\/} the range $m^3<N<(m+1)^3$, for what it is worth.
But, in the test case $N=m^3+m(2m+1)$ considered in \S 2, (22) gives
 the upper bound
$$ m+\frac{2m+1}{3m+1}. $$
Thus, (22) loses the advantage we found (1) has over (7) for such $N$
 (even if it remains better than (6)).

\section{Comparisons}

All comparisons, it is has often been said, are odious, but, as an
 anonymous reviewer wryly rejoined in the {\em Edinburgh Review\/}
 \cite[p.~400]{AR} for September, 1818:
\begin{quote}\small \em
No man, when he learns that the three angles of every triangle are
 equal to two right angles, ever thought of saying, that the series
 of comparisons by which that truth is demonstrated was invidious;
 neither has the fate of those interesting portions of space ever
 been deemed particularly hard, for having been subjected to such
 an investigation.
\end{quote}
The Greeks did debate the propriety of geometrical procedures --- we
 turn to one example in \S 5.3.
But their practical arithmetical competence was more pragmatic it
 seems.
Approximations tend to be stated blankly, without supporting argument,
 but also without comparison with other methods, as though truly, as
 Taisbak has it with (1), the Greeks {\em did not need any other
 corroboration than the fact that the method works\/}.

In contrast, for us today proposal of an approximative method is
 incomplete unless accompanied by examination of how well it performs
 against both rivals and the target.
So, in this section, we first look at an instance where Heron
 provides, not only a demonstration, but compares the resulting bound
 with an older rule of thumb; we then make a more thorough
 investigation of Enestr{\" o}m's identity; and we go on to show how
 a geometric scheme considered by Pappus can be adapted to improve
 on (1) for the family of numerical cases in \S 2.
We conclude by observing how the improving accuracy of (1), as
 revealed by (33), allows us to make good effect of rescaling
 (returns to scale).
The Newton-Raphson and Halley methods of approximating cube roots in
 (29) and (31), in contrast, do not guarantee such improving accuracy,
 even if some juggling may be possible (a rather more obvious
 distinction is that (1) is exact when $N$ is the cube of an integer).    

\subsection{{\em Metrica I.27--32}: Area of a circular segment}

Heron, in {\em Metrica I.27--32}, is concerned with formulae for
 the area of a circular segment (see \cite[II, pp.~330--331]{TLH}).
Let $AB$ be the arc of a circle subtending a segment less than a
 semicircle and let $C$ be the midpoint of the arc.
Then Heron asserts that the area subtended by $AB$ is greater than
 four thirds the area of the triangle $\triangle ABC$; that is,
 if the arc $AB$ has sagitta $h$ and subtended chord $b$, the
 subtended segment between arc and chord has area at least
\begin{equation}
 \frac{4}{3}\left(\frac{hb}{2}\right).
\end{equation}

But, rather out of character for him, Heron goes further, proving (23)
 in a manner reminiscent of Archimedes' {\em De quadratura parabolae,
 Prop.~24\/}.
However, despite being game to take on this task, Heron does not seem
 entirely sure of himself: he sets up his diagram as if intending to
 argue in one way, but then heads off in another; and underlying this
 dithering is a certain uneasiness in handling inequalities (at issue,
 in a sense, are returns to scale resulting from the circle's convexity,
 cf.\ \S 5.4).
So, it may be some surprise to find that, in {\em Metrica I.30, 31\/},
 Heron volunteers comparison of (23) with a more traditional
 approximation, namely
\begin{equation}
 \frac{h(b+h)}{2},
\end{equation}
 even stating, but without further comment, when one is to be
 preferred to the other.

This is all rather remarkable, and not unnaturally
 {\em Metrica I.27--32\/} has caught the attention of commentators.
Wilbur Knorr, in particular, has made much of the passage, returning to
 tease it out several times, as for example, in his books
 \cite[pp.~168--169]{WRK86} and \cite[pp.~498--501]{WRK89}, as well as
 in earlier papers on which the books build.
Knorr adjudicates the comparison of (23) and (24) in a footnote 
 \cite[p.~168, fn.~63]{WRK86} (in a further footnote
 \cite[p.~501, fn.~34]{WRK89}, he reports how advantage was not always
 taken of the improved bound):
\begin{quote}\em \small
[Hero] adds that one should use this rule when $b$ is less than
 three times $h$, but the former rule when $b$ is greater.
He does not explain this criterion, but one can see how it results
 from considering where the two rules yield the same result, namely,
$2bh/3=h(b+h)/2$, whence $b=3h$. \dots \\
The [former] rule, by virtue of its association with that for the
 parabolic segment, suggests an Archimedean origin.
One suspects that the rather sophisticated effort reported by Hero
 to assess the relative utility of these two rules for the circular
 segments is also due to an Archimedean insight.
\end{quote}

Now, there is no doubt that inequalities are more tricky to handle
 than equalities for pupils today, no less than in the past; and
 we all resort to simple means of reassurance that we have them
 right.
But, if Knorr's comments here arrest our attention, it is because
 of the incongruity between the supposed Archimedean origin of the
 comparison and the method advanced for seeing that it holds.
Perhaps Knorr is empathising too much with the difficulty Heron
 might have encountered in understanding some abstruse Archimedean
 proto-text.
Comparison of (23) and (24) would surely present little challenge
 to those, such as Archimedes, if not also Heron, for whom thinking
 in terms of areas was stock-in-trade.
  
In terms of areas, (23) tells us that the area of the subtended segment
 is a third more than the area of the triangle $\triangle ABC$, in
 keeping with the way the proof presented by Heron runs.
So, in place of (23), we might write the bound as
\begin{equation}
 \frac{hb}{2}+\frac{1}{3}\left(\frac{hb}{2}\right)=\frac{h(b+b/3)}{2}.
\end{equation}
Our areal intuition then suggests seeing in (24) and (25) triangles with
 common height $h$ and bases
$$ b+h;\quad b+\frac{b}{3}, $$
 respectively.
Which triangle has the larger area is simply a matter of which base is
 longer, leading to the conclusion that (25) is a better lower bound
 when the latter base is the larger, that is, when $b/3$ is bigger than
 $h$, as Heron claimed.

But, with Taisbak's stricture as our epigraph, the point to remember
 here --- and the point of this excursus --- is that this is only
 {\em our} intutition, not necessarily that of Heron or Archimedes,
 however plausible we fancy it to be.
On the other hand, they were clearly not in want of competence of their
 own.
 
\subsection{Enestr{\" o}m's identity}

It would be wrong to give the impression that the papers of Curtze
 \cite{MC} and Wertheim \cite{GW} are confined to the elaboration of
 Heron's text as discussed in the opening section.
For example, Curtze includes a list of quadratic approximations.
Wertheim anticipates the spirit of Taisbak in \cite{CMT}, providing a
 foundation on which Enestr{\" o}m builds in \cite{GHE2}.
Indeed, as Taisbak \cite{CMT2} playfully observes of any purported
 ``{\em new insight\/},'' on comparing Wertheim's contribution with
 his own,
\begin{quote}\em \small
If someone else said the same, it must be true.
If not, it is high time to have said it.
\end{quote}

Now, if we write
$$ \Delta_1=d_1-(n-m)^3;\quad \Delta_2=d_2-(m+1-n)^3, $$
 then Enestr{\" o}m, in \cite{GHE2}, goes through a series of
 algebraic manipulations that brings $n$ out in this notation as
\begin{equation}
 n=m+\frac{(m+1)\Delta_1}{(m+1)\Delta_1+m\Delta_2}.
\end{equation}
Clearly, if we ignore terms that are cubes of positive numbers
 less than unity, the right-hand side of (26) is just (1).
But (26) must hold as an identity, so going through a routine of
 solving for $n$, as Enestr{\" o}m does, might seem somewhat
 artificial.
Why not proceed more simply by direct computation with $\Delta_1$
 and $\Delta_2$?
We have
\begin{equation}
 \Delta_1=3mn(n-m);\quad \Delta_2=3(m+1)n(m+1-n),
\end{equation}
 expressions already familiar from \cite{CMT} as approximations
 for $d_1$ and $d_2$.
So, it readily follows that
$$ (m+1)^i\Delta_1+m^i\Delta_2=3m(m+1)n^i,\quad i=1, 2. $$ 
Hence (cf.\ (9), (11), (13) and (22)),
\begin{equation}
 n=\frac{(m+1)^2\Delta_1+m^2\Delta_2}{(m+1)\Delta_1+m\Delta_2}
  =m+\frac{(m+1)\Delta_1}{(m+1)\Delta_1+m\Delta_2},
\end{equation}
 as desired.

Looked at in this way, we see both that there is less mystery about
 Enestr{\" o}m's exact expression (26), but also less difference
 between him and later writers whose strategy is to get in early
 with the approximations for $d_1$ and $d_2$ given by (27), rather
 than waiting to the end.
Either way, while it is apparent that (1) {\em is\/} an approximation
 for the cube root of $N$, because we are modifying both numerator
 and denominator in the fraction we form in (28), we are left
 uncertain how {\em good\/} an approximation it is, or even whether
 we obtain an upper bound or a lower bound.
As Taisbak draws inspiration from the gradient of chords between
 successive integers and their cubes, his approach inherently sets
 up the expectation of an upper bound.

Naturally, a version of (28), and so of (1), can be developed for
 general intervals, as in \cite[p.~192]{WRK86} and \cite[p.~29, (1)]{DD}
 (that thoroughness is needed here can be seen from \cite[\S 2.1]{HLH}).
But Knorr's description in \cite[p.~192]{WRK86} of a prospective
 iterative application of such an extension of (1) also appears to
 be written in the expectation that the result gives an upper bound.
If, for some $a$ and $b$ not necessarily integers we have
 $a^3 < N < (a+b)^3$ and we obtain the approximation $a+b'$ after the
 manner of (1), as Knorr has us imagine, then certainly, at the next
 round of the iteration, we substitute for $a+b'$ for $a+b$, but only
 {\em if this approximation is an upper bound\/}.
In view of (21), we shall need to check this.
If, in the event, it turns out that $a+b'$ is a {\em lower\/} bound,
 we shall have to substitute it for $a$, {\em not} $a+b$, at the next
 round.

Knorr rightly goes on to question the authenticity of wiping away of
 small quantities, whenever in the scheme of things it happens,
 noting that we can reach the approximations in (27) in greater
 conformity with the Greek style by replacing the three terms on the
 left-hand side of (17) and (18) by three times their respective 
 middle terms, rather than being tied to versions of the binomial
 expansion (4) (see \cite[p.~193]{WRK86}).
So far as this approach goes, it is on a par with a Newton-Raphson
 approximation for the cube root of $N$, such as
\begin{equation}
 \frac{N+2m^3}{3m^2}
\end{equation}
 obtained by similarly replacing the same three terms by three times
 the last term, as Knorr also remarks.

\begin{figure}[htp]

\begin{center}

\includegraphics[width=15cm]{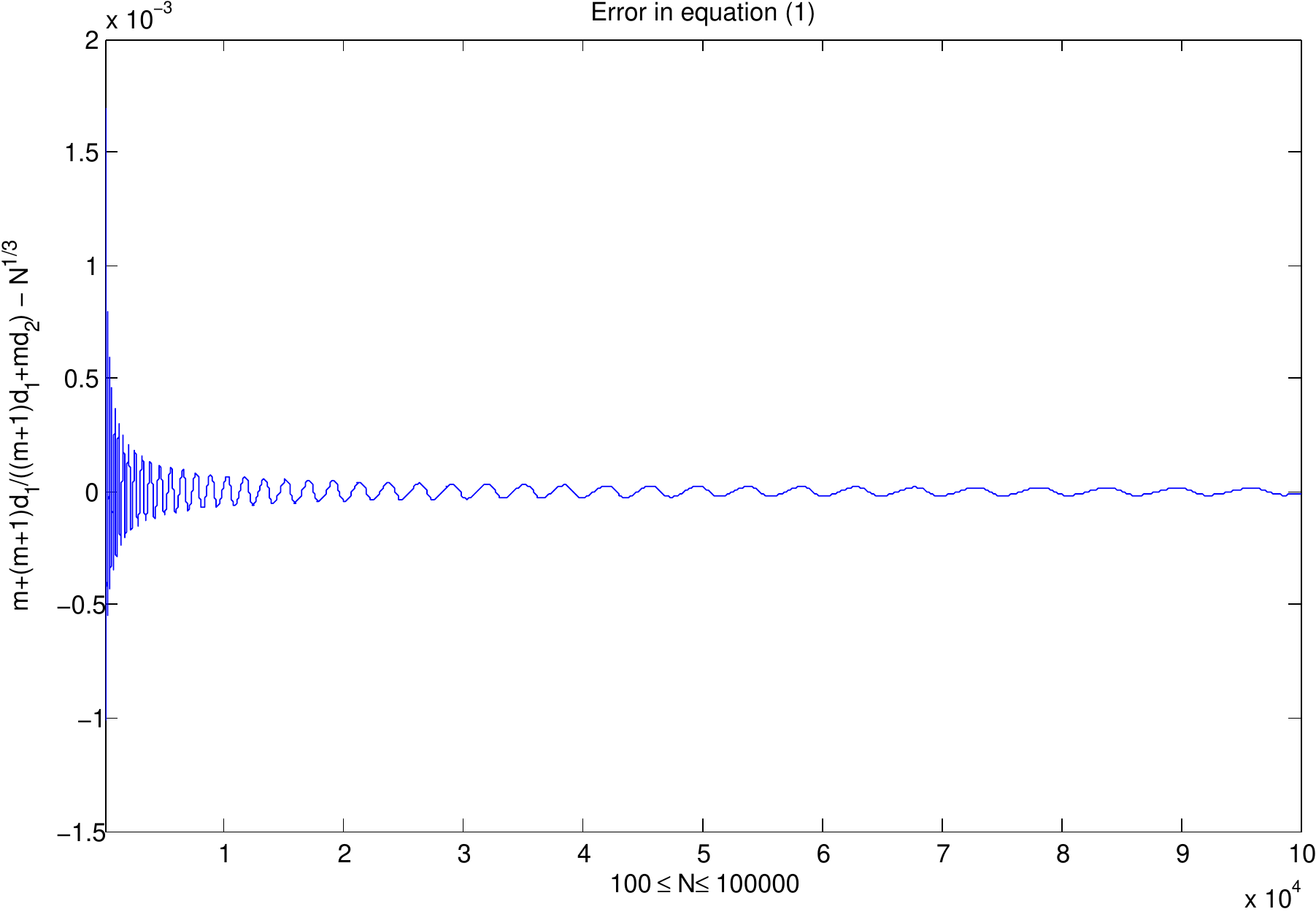}
\caption{Damped oscillation exhibited by error in (1), as given by (32)}

\end{center}
\end{figure}

For that matter, we could take this line of discussion further, by
 replacing the same three terms by three times the first term to
 obtain an approximation for the {\em square\/} of the cube root of $N$,
\begin{equation}
 \frac{2N+m^3}{3m},
\end{equation}
 and then cap this cleverness, by observing that an improved
 approximation for the cube root of $N$ proposed by Edmund Halley
 is given as the {\em ratio\/} of the expressions in (29) and (30):
\begin{equation}
  m\left(\frac{2N+m^3}{N+2m^3}\right).
\end{equation}
Halley's approximation in (31) does at least serve to remind us that
 in (1) we are also involved with a ratio, a ratio moreover, as (28)
 makes clear, of two blends of the approximations in (27).
Strangely enough, Knorr seems distracted from the significance of
 these differences between (1) and, say, (29), even while digressing
 at length on discoveries in approximation theory.

It may also be worth remembering that the statement of a result
 for illustrative purposes by way of a succinct algorithmic
 description, such as suits Heron's purpose in {\em Metrica III.20\/}
 might not be the formulation used were the result recast as a more
 formal proposition.
It is natural that historians of mathematics should wish to adhere to
 the text as they understand it, that is, to (1) as encapsulating the
 numerical instance in {\em Metrica III.20\/}; and that is what we
 find, with proposed proofs in which the manipulations of ratios
 closely follows the form of (1).
But, considering (13), (22), and now (28), in the general setting
 provided (8) and (9), we might suspect that it is these more symmetric
 equivalents of (1) that lend themselves more readily both to proof and
 to further examination.

\begin{figure}[htp]

\begin{center}

\includegraphics[width=15cm]{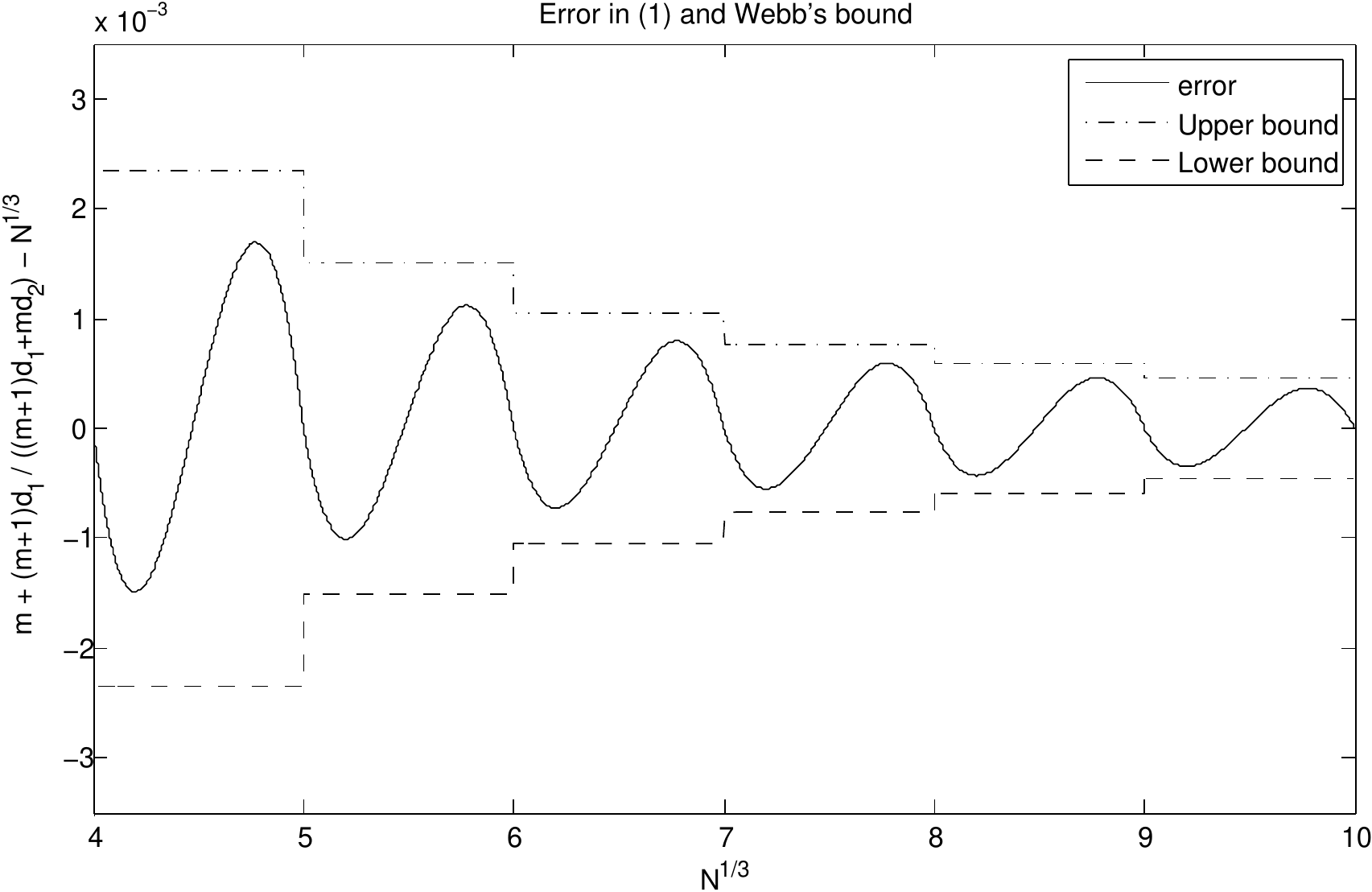}
\caption{Heron's Wave: error in (1) with Ward's bound superimposed}

\end{center}
\end{figure}

Thus, starting from (21), we find that
\begin{equation}
  \frac{(m+1)^{2}d_1+m^2d_2}{(m+1)d_1+md_2}-n=
 \frac{(n^2-m(m+1))(n-m)(m+1-n)}{(m+1)d_1+md_2}.
\end{equation}
To bound the absolute value of the left-hand side of (32) without
 going into too much fine detail, we note, first of all, that
$$ |n^2-m(m+1)| \leq m+1; $$
 secondly, by the inequality between geometric and arithmetic means
 (cf.\ {\em Elements VI.27\/})
$$ (n-m)(m+1-n)\leq \frac{1}{4}, $$
 with equality if and only if $n=m+1/2$; and thirdly
$$ (m+1)d_1+md_2 > m(d_1+d_2) \geq 3m^2(m+1). $$
Hence, putting these ingredients together, we conclude that
\begin{equation}
 \left|\frac{(m+1)^{2}d_1+m^2d_2}{(m+1)d_1+md_2}-n\right|<
  \frac{1}{12m^2},
\end{equation}
 of comparable order of magnitude to the bound $3/(80m^2)$ that
 Smyly tells us in \cite{JGS} had been obtained by Webb.
Another elementary bound is proved in \cite[Theorem 3]{DD}, but
 on the interval $(m,m+1)$ is is weaker than (33).

\subsection{{\em Synagogue III\/}: Two mean proportionals}  

Pappus musters in {\em Synagogue III\/} a collection of
 constructions of two mean proportionals between two line
 segments by non-planar means.
Perhaps by way of cautionary prologue, he also describes a
 geometrical solution, purportedly by plane considerations only,
 from some unnamed source, specifically with a view to showing
 that it fails.
The flaws in the construction are fairly transparent, and Pappus'
 demolition of them is not especially edifying.
However, for all the imperfections Pappus would have us see in
 it, the construction is not without other merits.
Knorr offers a sensitive geometrical re-appraisal at some length
 in \cite[pp. 64--70]{WRK89}; more recently, Serafina Cuomo has
 returned to the construction in a study \cite[\S 4.1]{SC} of
 Pappus' mathematics in the setting of Late Antiquity.
Earlier attempts at rehabilitating the construction tended to
 recast it as an iterative scheme of approximation to the mean
 proportionals, using an algebriac notation alien to the spirit
 of Pappus' {\em Synagogue\/}.
Nevertheless, what we might notice about this algebra for our
 present purposes is how well it meshes with the family of
 numerical examples in \S 2 generalising Heron's case, $N=100$,
 in {\em Metrica III.20}.

In this regard, the pioneering effort was made by Richard
 Pendlebury (1847--1902; Senior Wrangler, 1870) in a note
 \cite{RP} published in 1873, as reported in
 \cite[I, pp.~268--270]{TLH} (see further
 \cite[p.~64, fn.~8]{WRK89}; \cite[p.~130]{SC}). 
Suppose that $N=m^3-lm^2$, for some $l$ and $m$, then Pendlebury
 shows that iteration of the construction faulted by Pappus in
 {\em Synagogue III\/} can be generalised as a recursive
 computation,
\begin{equation}
  n_{i+1}=m-\frac{(m-n_i)lm^2}{m^3-n_i^3},
\end{equation}
 for some given $n_0$, with the $n_i$ successively better
 approximations to the cube root of $N$, giving upper bounds
 when $n_0$ is bigger than this cube root, and lower bounds
 when it is smaller.

Now, the family of $N$ in \S 2 generalising Heron's example is
 given by taking $l=1$.
If we start with our Heronian upper bound (5),
$$ n_0=m-\frac{m}{3m-1}=m(1-\frac{1}{3m-1}), $$
 then (34) gives the improved upper bound
\begin{equation}
  n_1=m-\frac{(3m-1)^2}{3(3m-1)(3m-2)+1}.
\end{equation}
In particular, for Heron's example, $N=100$ is the case $m=5$,
 when (35) yields
\begin{equation}
  n_1=5-\frac{196}{547}=4\frac{351}{547},
\end{equation} 
 an improvement on Heron's upper bound $4\frac{9}{14}$ for the
 cube root of 100.

In this exercise, we may be scrabbling after crumbs, waiting for a
 spark from heaven to fall.
This particular construction never seems to have attracted much
 attention until analysed by Pendlebury, although Leonardo Pisano
 and Gerolamo Cardano retained geometrical accounts of second mean
 proportionals in their discussions of cube root extraction.
But, over the course of countless Greek lives, there was presumably
 time for many other failed constructions and, in amongst them,
 some near-misses, possibly the occasional success --- after all,
 we still have Archimedes' {\em On the Measurement of a Circle\/}.

\subsection{Rescaling}

None of the ingredients we use in producing (33) could reasonably
 be said to be beyond the competence of the ancient Greek
 mathematicians, and yet we would naturally hesitate when it comes
 to an error bound like (33) itself.
Nevertheless, if we do have a sense that the going gets better,
 however we might come by it, we can always try rescaling.
Thus, to estimate the cube root in Heron's example, $N=100$, we
 might divide the estimate from (1) for the cube roots, say, of 800
 or 2700 by 2 or 3 respectively to get 
$$ 4\frac{322}{502}; \quad 4\frac{7328}{11421}; $$
 the first of these estimates is a lower bound not as close to
 the cube root of 100 as the upper bound in (36) while the second
 is an upper bound improving on that in (36).

Of course, (1) is most in error for some small values of $N$.
About the worst offender proportionately is $N=5$, when the estimate
 from (1) is $1\frac{8}{11}$, with a cube greater than 5.153.
It is here that we can use rescaling to good advantage.
Amusingly enough, if we divide the estimates from (1) for 40 or
 135 by 2 or 3 respectively, we come out with the same lower
 bound for the square root of 5, namely $1\frac{22}{31}$, with
 a cube greater than 4.997.
Going further and dividing the estimate from (1) for 320 by 4
 gives the upper bound $1\frac{615}{866}$, with a cube now less
 than 5.002.

Maybe there is some redemption to be found here, too, for the
 comparatively weak upper bound for square roots in (13), because,
 if we continue with the algebra there, we find that the
 diminution in the error is on the order of $1/m$.
For example, Heron, in {\em Metrica I.9\/}, wants to compute
 $\sqrt{1575}$ and notes he can get at this as
 $\frac{10}{2}\sqrt{63}$, offering the upper bound
 $7\frac{15}{16}$ for $\sqrt{63}$, either by Heronian iteration as
 in {\em Metrica I.8\/} or possibly as an application of (15)
 (cf.\ {\em Stereometrica I.33\/}).
Of course, if we stick with the same method and use it to approximate
 $\sqrt{1575}$ directly we come out with the same estimate either way.
However, as it so happens, Heron also alludes to $\sqrt{1575}$ in
 passing as ``{\em the square root of the fourth part of 6300\/}''
 (cf.\ \cite[p.~203]{EMB}).
But, if we divide the estimate of $\sqrt{6300}$ from (13) by 10,
 we obtain a (slightly) improved upper bound: $7\frac{1183}{1262}$.
Similarly, when Heron wants an approximation for $\sqrt{720}$ in
 {\em Metrica I.8\/}, his first estimate is the upper bound
 $26\frac{5}{6}$, whereas working (13) with $72,000$ improves this
 to $26\frac{30002}{36023}$.
 
Then, again, in any practical example, the convenience of working
 with an estimate may outweigh its accuracy, so such gains are
 largely a matter of theory.
Moreover, elsewhere, in {\em Geometrica 53, 54\/} 
 (cf. \cite[II, p.~321]{TLH}, when dealing with the 4-6-8 triangle,
 Heron seems to show some awareness that gains can be made from delay
 in the taking of square roots, initially proposing $a_1$, an upper
 bound with
$$ N = 4\sqrt{8\frac{7}{16}} < 11\frac{2}{3} = a_1, $$
 but then, on rewriting $N$ by multiplying into the square root,
 observing that we can do better using $a_2$, with
$$ N = \sqrt{135} < 11\frac{13}{21} = a_2. $$
Typically, nothing is said about the derivation of these bounds.
Interestingly enough though, Heronian iteration, as in (15),
 applied to $N$ gives $11\frac{5}{8}$, which falls in between the
 two bounds,
\begin{equation}
  a_1=11\frac{2}{3} > 11\frac{5}{8} > 11\frac{13}{21}=a_2;
\end{equation}
 $a_1$ results on applying Heronian iteration, or (15), to
 $\sqrt{136}=4\sqrt{8\frac{1}{2}}$; and $a_2$ improves on $a_1$
 precisely by Heronian iteration,
\begin{equation}
  a_2=\frac{1}{2}\left(\frac{135}{a_1}+a_1\right).
\end{equation}
A possible alternative derivation of $a_1$, in line with Heron's
 handing of $\sqrt{75}$ noted in \S 3.3, might be to stick with
 Heronian iteration in the form (14) for $N$, giving a less good
 upper bound $11\frac{14}{22}$, which, however, encourages nudging
 up to the simpler fraction $a_1$.
But all of this is speculative, and those who enjoy numerical
 coincidences will be amused to see the early Fibonacci numbers
 showing up in (37), still more perhaps to learn that these bounds
 are the 4th, 6th and 8th convergents of the continued fraction
 for $\sqrt{135}$.
Notice, however, that Heronian iteration with the middle bound in
 (37) yields $11\frac{307}{496}$, which does improve on $a_2$, if
 only just.

Thus, it is uncertain whether the improvement Heron notes here
 derives from {\em his\/} rescaling {\em per se\/} or from a
 change in the method of approximation.
Indeed, (38) may run slightly counter to the view in
 \cite[II, p.~326]{TLH} on Heron's own use of Heronian iteration,
 while leaving it a mystery as to how he obtained bounds that
 improve on a first instance of the method. 
Something similar might be at work in the handling of $\sqrt{28}$
 as discussed in \cite[p.~309]{EMB}.
In this case, we might expect the bound $5\frac{3}{10}$ (cf.\ (14)),
 but the weaker bound $5\frac{1}{3}$ (cf.\ (15)) lends itself more
 easily to improvement by Heronian iteration, giving $5\frac{7}{24}$.
However, what might require us to rethink, or at least re-express,
 the matter is the observation that rescaling combined with (15)
 does allow us to give the supposedly {\em improved\/} bounds in
 both cases more directly:
\begin{equation}
  \sqrt{28}=\frac{1}{3}\sqrt{252} <
 \frac{1}{3}\left(16-\frac{4}{32}\right)
 =\frac{1}{3}\left(15\frac{7}{8}\right)=5\frac{7}{24};
\end{equation}
\begin{equation}
  \sqrt{135}=\frac{1}{3}\sqrt{1215} < 
 \frac{1}{3}\left(35-\frac{10}{70}\right)
 =\frac{1}{3}\left(34\frac{6}{7}\right)=11\frac{13}{21}.
\end{equation}
Fortunately, under Taisbak's dispensation, we are not so pressed
 to account for the rather {\em weak\/} estimates Heron also uses
 on occasion, as, for example, $43\frac{1}{3}$ for $\sqrt{1875}$
 or $14\frac{1}{3}$ for $\sqrt{207}$, the former squaring to more
 than 1877, the latter to less than 206 (see
 \cite[II, pp.~326, 328]{TLH}).    
   
Smyly \cite[p.~67]{JGS}, in extolling the virtues of (1)
 for $N$ of the order of $10^6$ in comparison with tables of
 seven-figure logarithms, and Knorr \cite[p.~192]{WRK86}, in
 dilating on iterative use of (1), possibly overlook this simple
 trick of rescaling to obtain improved estimates for smaller $N$.
Scaling, in the elementary sense of the law of indices, is one
 thing; the notion of returns to scale another, rather more subtle. 
Some accounts of Greek approximations for $\sqrt{2}$ and $\sqrt{3}$
 would have us believe that the Greeks were great self-improvers,
 working their way to better estimates through solutions of Pell
 equations or the convergents of continued fractions which might
 be seen as implicitly involving a form of rescaling (indeed, not
 unlike (39) and (40)).
Taisbak \cite[\S 3]{CMT} asks in regard to his conjecture whether
 the Ancients knew and used sequences of differences.
With an eye to (39) and (40), we follow suit: did the Ancients know
 and use rescaling?

\section{A last reckoning}

Numerical corroboration, of course, might not be to everyone's taste.
Bartel Leendert van der Waerden (1903--1996), for one, in the original
 Dutch edition of {\em Ontwakende Wetenschap\/} ({\em Science
 Awakening\/}) \cite[p.~306]{vdW}, in 1950, places Heron in heavily
 weighted scales.
\begin{quote}\em \small
Laten we blij zijn, dat we de meesterwerken van Archimedes en
 Apollonios hebben, en niet treuren om het verlies van talloze
 rekenboekjes {\` a} la Heron.

[Let us rejoice in the masterworks of Archimedes and of Apollonius
 and not mourn the loss of numberless little accounting books after
 the manner of Heron.]
\end{quote}
The translation in English in 1954 is less pointed, but, recalling
 Heron's own mathematical outlook as expressed in the preface to
 {\em Metrica\/}, it is likely that he could at least hold his own
 (cf.\ \cite{JF}).

\bibliographystyle{plain}

\end{document}